\documentclass[12pt]{article}

\usepackage{amsmath}
\usepackage{theorem}
\usepackage{amssymb}
\usepackage{latexsym}
\usepackage{a4wide}

\newtheorem{thm}{Theorem}
\newtheorem{lem}[thm]{Lemma}
\newtheorem{prop}[thm]{Proposition}

{\theorembodyfont{\rmfamily} \newtheorem{exa}[thm]{Example}}
\newenvironment{rem}{\noindent{\bf Remark.}}{\newline}
\newenvironment{pf}{\noindent{\bf Proof.}}{\hbox{}\hfill $\Box$}

\newcommand{\Q}{\mathbb{Q}}
\newcommand{\R}{\mathbb{R}}
\newcommand{\Z}{\mathbb{Z}}
\newcommand{\mf}[1]{\mathfrak{#1}}
\newcommand{\hmf}[1]{\hat{\mathfrak{#1}}}
\newcommand{\emb}{\hookrightarrow}

\begin{document}

\title{Constructing homomorphisms between Verma modules}
\author{W. A. de Graaf, \\
School of Mathematics and Statistics,\\
University of Sydney, \\
Australia\\
email: {\tt wdg@maths.usyd.edu.au}}

\date{}
\maketitle

\begin{abstract}
We describe a practical method for constructing a nontrivial homomorphism 
between two Verma modules of an arbitrary semisimple Lie algebra. With some
additions the method generalises to the affine case.
\end{abstract}

A theorem of Verma, Bernstein-Gel'fand-Gel'fand gives a straightforward
criterion for the existence of a nontrivial homomorphism between Verma modules.
Moreover, the theorem states that such homomorphisms are always injective. In this
paper we consider the problem of explicitly constructing such a homomorphism 
if it exists. This boils down to constructing a certain element in the 
universal enveloping algebra of the negative part of the semismiple Lie algebra.\par
There are several methods known to solve this problem. Firstly, one can 
try and find explicit formulas. In this approach one fixes the type (but
not the rank). This has been carried out for type $A_n$ in \cite{mff},
Section 5, and for the similar problem in the quantum group case in \cite{dobrev3},
\cite{dobrev2}, \cite{dobrev1}.
In \cite{dobrev3} root systems of all types are considered, and the solution is given
relative to so-called straight roots, using a special basis of the universal
enveloping algebra (not of Poincar\'{e}-Birkhoff-Witt type). In \cite{dobrev2}, 
\cite{dobrev1} the solution
is given for types $A_n$ and $D_n$ for all roots, in a Poincar\'{e}-Birkhoff-Witt
basis. Our approach compares to this in the sense that we have an algorithm
that, given any root of a fixed root system, computes a general formula relative 
to any given Poincar\'{e}-Birkhoff-Witt basis (see Section \ref{sec3}).\par
A second approach is described in \cite{mff}, which gives a general construction of 
homomorphisms between Verma modules. However, it is not easy to see how to carry 
out this construction in practice. The method described here is a variant
of the construction in \cite{mff}, the difference being that we are able
to obtain the homomorphism explicitly.\par
In Section \ref{sec1} of this paper we review the theoretical
concepts and notation that we use, and describe the problem we deal with. In
Section \ref{sec2} we derive a few commutation formulas in the field of fractions
of $U(\mf{n}^-)$. Then in Section \ref{sec3} the construction of a homomorphism 
between Verma modules is described. In Section \ref{sec4} we briefly comment 
on the problem of finding compositions of inclusions. In Section \ref{sec5} we
comment on the analogous problem for affine algebras, and we show how our 
algorithm generalises to that case. Finally in Section \ref{sec6} we give an
application of the algorithm to the problem of constructing irreducible modules.
This is based on a result by P. Littelmann.\par
I have implemented the algorithms described in this paper in the computer algebra
system {\sf GAP}4 (\cite{gap4}). Sections \ref{sec3} and \ref{sec6} contain
tables of running times. All computations for these have been done on a PII 
600 Mhz processor, with 100M of memory of {\sf GAP}.

\section{Preliminaries}\label{sec1}

Let $\mf{g}$ be a semisimple Lie algebra, with root system $\Phi$, relative to a
Cartan subalgebra $\mf{h}$. We let $\Delta = \{\alpha_1,\ldots, \alpha_l\}$ be 
a fixed set of simple roots.
Let $\Phi^+ = \{ \alpha_1,\ldots ,\alpha_s\}$ be the set of positive roots (note that
here the simple roots are listed first).
Then there are root vectors $y_i=x_{-\alpha_i}$, $x_i=x_{\alpha_i}$ (for $1\leq i
\leq s$), and basis vectors $h_i \in \mf{h}$ (for $1\leq i\leq l$),  such that 
the set $\{x_1,\ldots, x_s,y_1,\ldots, y_s,h_1,\ldots, h_l\}$ forms a Chevalley 
basis of $\mf{g}$ (cf. \cite{hum}).
We have that $\mf{g} = \mf{n}^- \oplus \mf{h}\oplus \mf{n}^+$, where
$\mf{n}^-$, $\mf{n}^+$ are the subalgebras spanned by the $y_i$, $x_i$ 
respectively. \par
In the sequel, if $\beta=\alpha_i\in \Phi^+$, then we also write $y_{\beta}$ in place of
$y_i$.\par
We let $P$ denote the integral weight lattice spanned by the fundamental weights 
$\lambda_1,\ldots, \lambda_l$. Also $\Q P = \Q\lambda_1+\cdots +\Q \lambda_l$.
For $\lambda,\mu\in\Q P$ we write $\mu\leq \lambda$ if $\mu = \lambda - \sum_{i=1}^l
k_i \alpha_i$, where $k_i\in\Z_{\geq 0}$. Then $\leq $ is a partial order on
$\Q P$. \par
For $\alpha\in \Phi$ we have the reflection
$s_{\alpha} : \Q P\to \Q P$, given by $s_{\alpha}(\lambda) = \lambda -\langle \lambda,
\alpha^{\vee}\rangle \alpha$.  \par
Let $U(\mf{g})$ denote the universal enveloping algebra of $\mf{g}$. We consider 
$U(\mf{g})$ as a $\mf{g}$-module by left multiplication. Let $\lambda = \sum a_i
\lambda_i \in \Q P$, and let $J(\lambda)$ be the $\mf{g}$-submodule of $U(\mf{g})$ 
generated by  $h_i - a_i+1$ for $1\leq i\leq l$ and $x_i$ for $1\leq i\leq s$. 
Then $M(\lambda) = U(\mf{g})/J(\lambda)$ is a $\mf{g}$-module. It is called a Verma
module. As $U(\mf{g}) = 
U(\mf{n}^-)\oplus J(\lambda)$ we see that $U(\mf{n}^-)\cong M(\lambda)$ (as
$U(\mf{n}^-)$-modules). Let $v_{\lambda}$ denote the image of $1$ under this 
isomorphism. Then $h_i\cdot v_{\lambda} = (a_i-1) v_{\lambda}$, and $x_i\cdot 
v_{\lambda}=0$. Furthermore, all other elements of $M(\lambda)$ can be written
as $Y\cdot v_{\lambda}$, where $Y\in U(\mf{n}^-)$. \par
Let $\nu = \sum_{i=1}^l k_i\alpha_i$, where $k_i\in \Z_{\geq 0}$. Then we let
$U(\mf{n}^-)_{\nu}$ be the span of all $y_{i_1}\cdots y_{i_r}$ such that 
$\alpha_{i_1}+\cdots +\alpha_{i_r} = \nu$. \par
For a proof of the following theorem we refer to \cite{bgg}, \cite{dix}.

\begin{thm}[Verma, Bernstein-Gel'fand-Gel'fand]\label{thm1}
Let $\lambda,\mu\in \Q P$, and set 
$$R_{\mu,\lambda}= {\rm Hom}_{U(\mf{g})}(M(\mu),M(\lambda)).$$
Then
\begin{enumerate}
\item $\dim R_{\mu,\lambda}\leq 1$,
\item non-trivial elements of $R_{\mu,\lambda}$ are injective,
\item $\dim R_{\mu,\lambda}=1$ if and only if there are positive roots 
$\alpha_{i_1},\ldots, \alpha_{i_k}$ such that 
$$\mu \leq s_{\alpha_{i_1}}(\mu)\leq s_{\alpha_{i_2}}s_{\alpha_{i_1}}(\mu)\leq
\cdots \leq s_{\alpha_{i_k}}\cdots s_{\alpha_{i_1}}(\mu)=\lambda.$$
\end{enumerate}
\end{thm}

The problem we consider is to construct a non-trivial element in $R_{\mu,\lambda}$ if
$\dim R_{\mu,\lambda}=1$. By Theorem \ref{thm1}, this boils down to finding an element
in $R_{\mu,\lambda}$ if $\mu = s_{\alpha}(\lambda)=\lambda - \langle \lambda, 
\alpha^{\vee}\rangle
\alpha$ and $\langle \lambda, \alpha^{\vee}\rangle \in \Z_{> 0}$. Suppose that
we are in this situation, and set $h=\langle \lambda, \alpha^{\vee}\rangle$.
An element $Y\cdot v_{\lambda}\in M(\lambda)$, where $Y\in U(\mf{n}^-)$
is said to be {\em singular} if
$x_{\alpha}\cdot (Y\cdot v_{\lambda})=0$ for $\alpha\in\Phi^+$.
Let $\psi \in R_{\mu,\lambda}$ be a non-trivial $U(\mf{g})$-homomorphism.
Then $\psi(v_{\mu})=Y\cdot v_{\lambda}$ for some $Y\in U(\mf{n^-})$ with
$Y\cdot v_{\lambda}$ singular. We have $h_i y = yh_i -\langle\nu, \alpha_i^{\vee}
\rangle y$ for all $y\in U(\mf{n}^-)_{\nu}$. Hence $h_i\cdot (y\cdot v_{\lambda}) =
(\langle \lambda-\nu, \alpha_i^{\vee}\rangle -1)yv_{\lambda}$. So, as 
$h_i\cdot (Yv_{\lambda}) = (\langle \mu, \alpha_i^{\vee}\rangle -1)Yv_{\lambda}$
we see that $Y\in U(\mf{n}^-)_{h\alpha}$.
Conversely, if we have a 
$Y\in U(\mf{n^-})_{h\alpha}$ such that $Y\cdot v_{\lambda}$ is singular, 
then $\psi : M(\mu)\to M(\lambda)$ defined by 
$\psi( Y'\cdot v_{\mu} ) = Y'Y\cdot v_{\lambda}$ will be a non-trivial element
of $R_{\mu,\lambda}$. So the problem reduces to finding a $Y\in U(\mf{n^-})_{h\alpha}$ 
such that
$Y\cdot v_{\lambda}$ is singular. Note that this can be done by writing down a basis
for $U(\mf{n^-})_{h\alpha}$ and computing a set of linear equations for $Y$. However, this
algorithm becomes rather cumbersome if $\dim U(\mf{n^-})_{h\alpha}$ gets large. We will
describe a more direct method.

\section{The field of fractions}\label{sec2}

From \cite{dix}, \S 3.6 we recall that $U(\mf{n}^-)$ has a (non-commutative) 
field of fractions, denoted by $K(\mf{n}^-)$. It consists of all elements
$ab^{-1}$ for $a\in U(\mf{n}^-)$, $b\in U(\mf{n}^-)\setminus\{0\}$. For the 
definitions of addition and multiplication in $K(\mf{n}^-)$ we refer to \cite{dix},
\S 3.6. They imply $aa^{-1}=a^{-1}a=1$. \par
Let $\alpha,\beta\in\Phi^+$. If $\alpha+\beta\in\Phi^+$ then we let $N_{\alpha,\beta}$
be the scalar such that $[y_{\alpha},y_{\beta}] = -N_{\alpha,\beta} y_{\alpha+\beta}$.
Also set $P_{\alpha,\beta}= \{i\alpha +j\beta\mid i,j\geq 0\} \cap \Phi^+$.
Then there are seven possibilities for $P_{\alpha,\beta}$:
\begin{itemize}
\item[(I)] $P_{\alpha,\beta} = \{\alpha,\beta\}$,
\item[(II)] $P_{\alpha,\beta} = \{\alpha,\beta,\alpha+\beta\}$
\item[(III)] $P_{\alpha,\beta} = \{\alpha,\beta,\alpha+\beta,\alpha+2\beta\}$,
\item[(IV)] $P_{\alpha,\beta} = \{\alpha,\beta,\alpha+\beta,2\alpha+\beta\}$,
\item[(V)] $P_{\alpha,\beta} = \{\alpha,\beta,\alpha+\beta,2\alpha+\beta,
3\alpha+\beta,3\alpha+2\beta\}$
\item[(VI)] $P_{\alpha,\beta} = \{\alpha,\beta,\alpha+\beta,\alpha+2\beta,
\alpha+3\beta,2\alpha+3\beta\}$,
\item[(VII)] $P_{\alpha,\beta} = \{\alpha,\beta,\alpha+\beta,2\alpha+\beta,
\alpha+2\beta\}$.
\end{itemize}

\begin{lem}\label{lem1}
In case (I) we have $y_{\beta}^my_{\alpha}^n = y_{\alpha}^n y_{\beta}^m$ for all
$m,n\in\Z$.
\end{lem}

\begin{pf}
If $n>0$ then $y_{\beta}y_{\alpha}^n = y_{\alpha}^n y_{\beta}$. Multiplying this
relation on the left and on the right by $y_{\alpha}^{-n}$ we get 
$y_{\beta}y_{\alpha}^{-n} = y_{\alpha}^{-n} y_{\beta}$. So we have $y_{\beta}
y_{\alpha}^n = y_{\alpha}^n y_{\beta}$ for all $n\in \Z$. From this it follows that
$y_{\beta}^my_{\alpha}^n = y_{\alpha}^n y_{\beta}^m$ for $m>0$, $n\in \Z$. If we now
multiply this from the left and the right by $y_{\beta}^{-m}$ we get the result
for $m<0$ as well.
\end{pf}

Since
$$\binom{n}{k} = \frac{n(n-1)\cdots (n-k+1)}{k!},$$
these binomial coefficients are defined for arbitrary $n\in \Q$, and $k\in 
\Z_{\geq 0}$. In fact, we see that $\binom{n}{k}$ is a polynomial of degree 
$k$ in $n$. Note also that if $n\in \Z$ and $0\leq n <k$ then the coefficient
is $0$.

\begin{lem}\label{lem2}
In case (II) we have for $m\geq 0$, $n\in\Z$,
$$y_{\beta}^my_{\alpha}^n = \sum_{k=0}^m N_{\alpha,\beta}^k 
\binom{m}{k}\binom{n}{k} k! y_{\alpha}^{n-k}y_{\beta}^{m-k}y_{\alpha+\beta}^k.$$
\end{lem}

\begin{pf}
First of all, this formula is known for $m,n\geq 0$ (see, e.g., \cite{gra5}).
In particular, for $n>0$ we have $y_{\beta} y_{\alpha}^n = y_{\alpha}^ny_{\beta}
+ N_{\alpha,\beta}n y_{\alpha}^{n-1}y_{\alpha+\beta}$. If we multiply this 
relation on the left and the right with $y_{\alpha}^{-n}$, and use Lemma
\ref{lem1}, then we get 
it for all $n\in \Z$. Now the formula for $m>1$ is proved by induction. 
\end{pf}

\begin{lem}\label{lem3}
In case (III) we have for $m\geq 0$, $n\in\Z$,
$$y_{\beta}^my_{\alpha}^n = \sum_{\substack{ k,l\geq 0\\k+2l\leq m}} c_{k,l}^{m,n}
y_{\alpha}^{n-k-l} y_{\beta}^{m-k-2l} y_{\alpha+\beta}^k y_{\alpha+2\beta}^l,$$
where
$$c_{k,l}^{m,n} = N_{\alpha,\beta}^{k+l} (\frac{1}{2}N_{\beta,\alpha+\beta})^l 
\binom{n}{k+l}\binom{m}{k+2l}\binom{k+l}{l} (k+2l)!.$$
\end{lem}

\begin{pf}
This goes in exactly the same way as the proof of Lemma \ref{lem2}.
\end{pf}

\begin{lem}\label{lem4}
In case (IV) we have for $m\geq 0$, $n\in\Z$,
$$y_{\beta}^my_{\alpha}^n = \sum_{\substack{ k,l\geq 0\\k+l\leq m}} c_{k,l}^{m,n}
y_{\alpha}^{n-k-2l} y_{\beta}^{m-k-l} y_{\alpha+\beta}^k y_{2\alpha+\beta}^l,$$
where 
$$c_{k,l}^{m,n} =N_{\alpha,\beta}^{k+l} (\frac{1}{2}N_{\alpha,\alpha+\beta})^l 
\binom{n}{k+2l}\binom{m}{k+l}\binom{k+l}{l} (k+2l)!.$$
\end{lem}

\begin{pf}
Again we get the formula for $m,n\geq 0$ from \cite{gra5}. In this case the 
formula for $m=1$, $n\geq 0$ reads
$$y_{\beta}y_{\alpha}^n = y_{\alpha}^ny_{\beta} + N_{\alpha,\beta} n y_{\alpha}^{n-1}
y_{\alpha+\beta} + N_{\alpha,\beta} N_{\alpha,\alpha+\beta} \binom{n}{2} 
y_{\alpha}^{n-2} y_{2\alpha+\beta}.$$
If we multiply this on the left and the right by $y_{\alpha}^{-n}$, and use 
Lemmas \ref{lem1}, \ref{lem2} we get the same relation with $n$ replaced by $-n$.
So the case $m=1$, $n\in \Z$ follows. The formula for $m>1$ now follows by induction.   
\end{pf}

The cases (V), (VI), (VII) can only occur when the root system has a component 
of type $G_2$. We omit the formulas for these cases; they can easily be derived 
from those given in \cite{gra5}. \par
Now let $a=y_1^{n_1}\cdots y_s^{n_s}$ be a monomial in 
$U(\mf{n^-})$. For $\beta\in\Phi^+$ and $m,n\in \Z$ consider the element
$y_{\beta}^may_{\beta}^{-n}$. By repeatedly applying Lemmas \ref{lem1},
\ref{lem2}, \ref{lem3}, \ref{lem4} we see that
\begin{equation}\label{eq2.1}
y_{\beta}^{m}ay_{\beta}^{-n} = \sum_{(k_1,\ldots, k_t)\in I}
c(k_1,\ldots,k_t) y_{\beta}^{m-n-p_1k_1-\cdots -p_tk_t} a(k_1,\ldots, k_t).
\end{equation}
Here the $a(k_1,\ldots, k_t)\in U(\mf{n}^-)$, 
the (finite) index set $I$, the $p_i\in \Z_{>0}$ are all independent of $n$, they
only depend on $a$. Only the exponents of $y_{\beta}$  and the coefficients 
$c(k_1,\ldots,k_t)$ (which are polynomials in $n$) depend on $n$. \par
Now we take $m,n\in \Q$ such that $m-n\in \Z$. Then we define 
$y_{\beta}^{m}ay_{\beta}^{-n}$
to be the right-hand side of (\ref{eq2.1}), and we say that $y_{\beta}^{m}ay_{\beta}^{-n}$
is an element of $K(\mf{n^-})$. More generally, if $Y$ is a linear combination
of monomials, and $m,n\in\Q$ such that $m-n\in \Z$ then 
$y_{\beta}^{m}Yy_{\beta}^{-n}$ is an element of $K(\mf{n^-})$.

\section{Constructing singular vectors}\label{sec3}

Here we suppose that we are given a $\lambda\in\Q P$ and $\alpha\in \Phi^+$ with
$\langle \lambda,\alpha^{\vee}\rangle =h\in\Z_{>0}$. The problem is to find a
$Y\in U(\mf{n^-})_{h\alpha}$ such that $Y\cdot v_{\lambda}$ is a singular vector.\par
We recall that $l=|\Delta|$ is the rank of the root system. Let $1\leq i\leq l$, then
\begin{equation}\label{eq3.1}
x_iy_i^r\cdot v_{\lambda} = r(\langle \lambda, \alpha_i^{\vee}\rangle -r)y_i^{r-1}\cdot
v_{\lambda}.
\end{equation}

\begin{lem}\label{lem3.1}
Suppose that $\alpha\in\Delta$, i.e., $\alpha=\alpha_i$, $1\leq i\leq l$.
Then $y_i^h\cdot v_{\lambda}$ is a singular vector. 
\end{lem}

\begin{pf}
This follows from (\ref{eq3.1}), cf. the proof of \cite{bgg}, Lemma 2.
\end{pf}

Note that this solves the problem when $\mf{g}=\mf{sl}_2$. So in the remainder 
we will assume that the rank of the root system is at least $2$. By an
embedding $\phi : M(\mu)\emb M(\lambda)$ we will always mean an
injective $U(\mf{g})$-homomorphism.

\begin{lem}\label{lem3.2}
Suppose that $\nu,\eta\in P$, and $\beta\in\Delta$ is such that 
$m=\langle\nu,\beta^{\vee}\rangle$ is a non-negative integer. Suppose further that
we have an embedding $\psi : M(\nu)\emb M(\eta)$ given by $\psi(v_{\nu})=
Yv_{\eta}$. Set $n=\langle \eta, \beta^{\vee}\rangle$. Then
$y_{\beta}^m Y y_{\beta}^{-n}$ is an element of $U(\mf{n}^-)$ and
we have an embedding $\phi : M(s_{\beta}\nu)\emb M(s_{\beta}\eta)$ given by
$\phi( v_{s_{\beta}\nu} ) =  y_{\beta}^m Y y_{\beta}^{-n} \cdot v_{s_{\beta}\eta}$.
\end{lem}

\begin{pf}
If $n\leq 0$ then the first statement is clear. The embedding $\phi$ is the 
composition $M(s_{\beta}\nu)\emb M(\nu)\emb M(\eta)\emb M(s_{\beta}\eta)$,
where the first and the third maps follow from Lemma \ref{lem3.1}.\par
If $n>0$, then we view $M(s_{\beta}\eta)$ as a submodule of $M(\eta)$. We have
$v_{s_{\beta}\eta}=y_{\beta}^nv_{\eta}$ (Lemma \ref{lem3.1}). Set $v=y_{\beta}^m Y
v_{\eta}$; then $v$ is a singular vector (being the image of $v_{s_{\beta}\nu}$ under
$M(s_{\beta}\nu)\emb M(\nu)\emb M(\eta)$). We claim that $v\in M(s_{\beta}\eta)$.
Suppose that this claim is proved. Then there is a $Y'\in U(\mf{n}^-)$ such that
$v=Y'v_{s_{\beta}\eta}$. But that means that $y_{\beta}^m Y = Y'y_{\beta}^n$, and
the lemma follows.\par
The claim above is proved in \cite{bgg}. For the sake of completeness we transcribe
the argument. Set $V=M(\eta)/M(s_{\beta}\eta)$, and let
$\bar{v}_{\nu}$ denote the image of $\psi(v_{\nu})$ in $V$; then $\bar{v}_{\nu} =
X\cdot \bar{v}_{\eta}$, for some $X\in U(\mf{n}^-)$.
For $k\geq 0$ write $y_{\beta}^kX=X_1y_{\beta}^{k_1}$.
By increasing $k$ we can get $k_1$ arbitrarily large (cf. \cite{dix}, Lemma 7.6.9;
it also follows by straightforward weight considerations). By Lemma \ref{lem3.1}
we know that $y_{\beta}^n v_{\eta}\in M(s_{\beta}\eta)$. Therefore there is a
$k>0$ such that $y_{\beta}^k\bar{v}_{\nu}=0$. Then by using (\ref{eq3.1}) we see that 
the smallest such $k$ must be equal to $m$.
\end{pf}

\begin{prop}\label{prop3.1}
Let $\nu,\eta\in\Q P$ be such that $\nu = s_{\gamma}(\eta)=\eta-k\gamma$, where
$\gamma\in\Phi^+$ and $k\in\Z_{>0}$. 
Let $Y\in U(\mf{n}^-)_{k\gamma}$ be such that $Y\cdot v_{\eta}$ is singular.
Let $\beta\in\Delta$, $\beta\neq\gamma$. 
Set $m=\langle \nu, \beta^{\vee}\rangle$, $n=\langle \eta, \beta^{\vee}\rangle$.
Then $y_{\beta}^mYy_{\beta}^{-n}$ is an element of $K(\mf{n}^-)$; it is even an
element of $U(\mf{n}^-)$. Secondly, we have an embedding $\phi : M(s_{\beta}\nu)\emb 
M(s_{\beta}\eta)$ given by $\phi( v_{s_{\beta}\nu} ) =  y_{\beta}^m Y y_{\beta}^{-n} 
\cdot v_{s_{\beta}\eta}$.
\end{prop}

\begin{pf}
We have that $m-n = -k\langle \gamma,\beta^{\vee}\rangle\in\Z$, so
$y_{\beta}^mYy_{\beta}^{-n}$ is an element of $K(\mf{n}^-)$. \par
Set $V=\{ \mu\in\Q P\mid \langle\mu, \gamma^{\vee}\rangle=k\}$, which is a 
hyperplane in $\Q P$, containing $\eta$. 
Let $\{a_1,\ldots, a_t\}$ be a basis of $U(\mf{n}^-)_{k\gamma}$.
Take $\mu = \sum_{i=1}^l r_i \lambda_i \in V$ and set $\tilde{\mu} = s_{\gamma}(\mu)
=\mu-k\gamma$. Then by Theorem \ref{thm1} there is a $Y_{\mu}=\sum_{i=1}^t \zeta_i a_i$
such that $Y_{\mu}\cdot v_{\mu}$ is singular. Here the $\zeta_i$ are polynomial 
functions of the $r_i$. (Indeed, if we write linear equations for the $Y_{\mu}$, then
the coefficients depend linearly on the $r_i$. Hence the coefficients of a solution
are polynomial functions of the $r_i$.) \par
Set $p=\langle \tilde{\mu},\beta^{\vee}\rangle$, $q=\langle \mu,\beta^{\vee}\rangle$.
Then $Y' =y_{\beta}^p Y_{\mu} y_{\beta}^{-q} =\sum_j c_j b_j$, where the $b_j$
are linearly independent elements of $K(\mf{n}^-)$, and the $c_j$ are coefficients
that depend polynomially on the $r_i$. Now Lemma \ref{lem3.2} implies that
if the $r_i\in\Z$ and $p\geq 0$, then $Y'\in U(\mf{n}^-)$.  Let now $j$ be such that
$b_j\not\in U(\mf{n}^-)$. If  the $r_i\in\Z$ and $p\geq 0$, then $c_j=0$. 
Suppose that $\beta=\alpha_{i_0}$, the $i_0$-th simple root. Then the 
condition $p\geq 0$ amounts
to $r_{i_0} \geq k\langle \gamma,\beta^{\vee}\rangle$. We have that $\mu\in V$ if and
only if $\sum_{i=1}^l u_i r_i=k$, where the $u_i$ are certain elements of $\Z$. 
Also, since $\beta\neq \gamma$ at least one $u_i\neq 0$ with $i\neq i_0$. We see that
the requirement $r_{i_0} \geq k\langle \gamma,\beta^{\vee}\rangle$ cuts a half space $W$ 
off $V$. Furthermore $V\cap P$ is an $(l-1)$-dimensional lattice in $V$ (cf. \cite{bgg}).
The conclusion is that $c_j=0$ if $\mu \in W\cap P$. Since the $c_j$ are polynomials 
in the $r_i$, it follows that $c_j=0$ if $\mu\in V$. In particular, $y_{\beta}^m Y
y_{\beta}^{-n}$ lies in $U(\mf{n}^-)$.\par
Finally we note that $Y'\cdot v_{s_{\beta}\mu}$ is singular, by the same arguments.
(Indeed, $x_i\cdot (Y'\cdot v_{s_{\beta}\mu}) = \sum_j f_j z_j\cdot v_{s_{\beta}\mu}$
where the $f_j$ are polynomials in the $r_i$, and the $z_j$ are elements of 
$U(\mf{n}^-)$. Since the $f_j$ are zero when $\mu\in W\cap P$ we have that $f_j=0$
when $\mu\in V$.) In particular, $y_{\beta}^m Y y_{\beta}^{-n}\cdot v_{s_{\beta}
\eta}$ is singular.
\end{pf}

\begin{exa}\label{exa1}
To illustrate the argument in the preceding proof,
consider the Lie algebra of type $A_3$, with simple roots $\alpha,\beta,\gamma$
(with $\beta$ corresponding to the middle node of the Dynkin diagram). 
Then it is possible to choose a Chevalley basis such that
$[y_{\alpha},y_{\beta}] = y_{\alpha+\beta}$, $[y_{\alpha},y_{\beta+\gamma}]
=y_{\alpha+\beta+\gamma}$, $[y_{\beta},y_{\gamma}]=y_{\beta+\gamma}$, 
$[y_{\gamma},y_{\alpha+\beta}]= - y_{\alpha+\beta+\gamma}$.
Set $a_1=y_{\alpha}y_{\beta}y_{\gamma}$, $a_2=y_{\gamma}y_{\alpha+\beta}$,
$a_3=y_{\alpha}y_{\beta+\gamma}$, $a_4=y_{\alpha+\beta+\gamma}$. Then
$\{a_1,a_2,a_3,a_4\}$ is a basis of $U(\mf{n}^-)_{\alpha+\beta+\gamma}$.\par
We abbreviate a weight $r_1\lambda_1+r_2\lambda_2+r_3\lambda_3$ by $(r_1,r_2,r_3)$. 
Let $V$ be the hyperplane in $\Q P$ consisting of all weights $\mu$ such that 
$\langle \mu,(\alpha+\beta+\gamma)^{\vee}\rangle =1 $, i.e., $V=\{
(r_1,r_2,r_3) \mid r_1+r_2+r_3=1\}$. Let $\mu=(r_1,r_2,r_3)\in V$ and set 
$\tilde{\mu} = s_{\alpha+\beta+\gamma}(\mu) = (r_1-1,r_2,r_3-1)$. Set
$Y_{\mu} = a_1-r_1a_2-(r_1+r_2)a_3-r_1r_3a_4$; then $Y_{\mu}\cdot v_{\mu}$ is
singular. Set $p=\langle \tilde{\mu},\alpha^{\vee}\rangle=r_1-1$ and 
$q=\langle \mu,\alpha^{\vee}\rangle=r_1$. Now $Y' = y_{\alpha}^p Y_{\mu} 
y_{\alpha}^{-q} = y_{\beta}y_{\gamma}-(r_1+r_2)y_{\beta+\gamma} +
r_1(1-r_1-r_2-r_3)y_{\alpha}^{-1} y_{\alpha+\beta+\gamma}$. According
to Lemma \ref{lem3.2} this is an element of $U(\mf{n}^-)$ whenever
$(r_1,r_2,r_3)\in V$ with the $r_i$ integral and $p\geq 0$. Therefore the 
coefficient of $y_{\alpha}^{-1} y_{\alpha+\beta+\gamma}$ has to vanish,
which is indeed the case. We see that $Y'$ lies in $U(\mf{n}^-)$ for all
$(r_1,r_2,r_3)\in V$.
\end{exa}

Now we return to the situation of the beginning of the section. We have $\lambda\in
\Q P$, $\alpha\in\Phi^+$ with $\langle \lambda,\alpha^{\vee}\rangle = h \in \Z_{>0}$.
Set $\mu = s_{\alpha}(\lambda)=\lambda-h\alpha$. 
To obtain an embedding $M(\mu)\emb M(\lambda)$, we perform the following steps:
\begin{enumerate}
\item Select $\beta_1,\ldots, \beta_r\in\Delta$ and positive roots 
$\alpha_0,\ldots, \alpha_r$ in the following way. Set $\alpha_0=\alpha$, and 
$k=0$. Then:
\begin{enumerate}
\item If $\alpha_k\in \Delta$, then set $r=k$ and go to step 2.
\item Otherwise, let $\beta_{k+1}\in\Delta$ be such that $\langle \alpha_k, \beta_{k+1}
^{\vee}\rangle >0$, and set $\alpha_{k+1} = s_{\beta_{k+1}}(\alpha_k)$, and 
$k:= k+1$. Return to (a).
\end{enumerate}
\item Set $\beta=\alpha_r\in\Delta$. 
 For $1\leq k\leq r$ set $a_k = -\langle \mu, s_{\beta_1}\cdots s_{\beta_{k-1}}
(\beta_k)^{\vee}\rangle$, and $b_k = \langle \lambda, s_{\beta_1}\cdots s_{\beta_{k-1}}
(\beta_k)^{\vee}\rangle$. 
\item Set $Y_0 = y_{\beta}^h$, and for $0\leq k\leq r-1$:
$$Y_{k+1} = y_{\beta_{r-k}}^{a_{r-k}} Y_{k} y_{\beta_{r-k}}^{b_{r-k}}.$$
\end{enumerate}

\begin{rem}
Note that the $\beta_{k+1}$ in step 1 (b) exists because otherwise $\langle \alpha_k,
\gamma^{\vee}\rangle \leq 0$ for all $\gamma\in \Delta$, and this implies
that the set $\Delta \cup \{\alpha_k\}$ is linearly independent
(cf. \cite{jac}, Chapter IV, Lemma 1), which is not possible
since $\alpha_k\not\in\Delta$. Also, all $\alpha_k$ must be positive roots
because $s_{\gamma}$ permutes the positive roots other than $\gamma$, for $\gamma
\in\Delta$. Then the loop in 1. must terminate because the height of $\alpha_k$
decreases every step. 
\end{rem}

\begin{prop}\label{prop3.2}
All $Y_k$ are elements of $U(\mf{n}^-)$ and we have an embedding $M(\mu)\emb M(\lambda)$
given by $v_{\mu}\mapsto Y_r\cdot v_{\lambda}$.
\end{prop}

\begin{pf}
We write $s_i=s_{\beta_i}$. For $0\leq k\leq r$ we set $w_k = s_{r-k}\cdots s_1$
(so $w_r=1$), and 
$\mu_k = w_k\mu$, $\lambda_k = w_k\lambda$. We claim that there is an embedding
$M(\mu_k)\emb M(\lambda_k)$ given by $v_{\mu_k}\mapsto Y_k \cdot v_{\lambda_k}$. First
we look at the case $k=0$. Note that $s_r\cdots s_1 (\alpha)=\beta\in\Delta$.
Since for $w$ in the Weyl group we have $ws_{\beta}w^{-1} = s_{w\beta}$
we get $s_{\alpha}=s_1\cdots s_rs_{\beta}s_r\cdots s_1=w_0^{-1}s_{\beta}w_0$. 
Therefore $\mu_0 = w_0s_{\alpha}(\lambda) = s_{\beta}
(\lambda_0)$, and $\langle \lambda_0, \beta^{\vee}\rangle = \langle \lambda,
s_1\cdots s_r(\beta)^{\vee}\rangle = \langle \lambda, \alpha^{\vee}\rangle =h$. The
case $k=0$ now follows by Lemma \ref{lem3.1}.\par
Now suppose we have an embedding $M(\mu_k)\emb M(\lambda_k)$ as above. Note that
$w_{k+1}=s_{r-k}w_k$ and $\alpha_k=w_k\alpha$.
Also $\mu_k=w_k\mu = \lambda_k -h \alpha_k$, and
$\langle \lambda_k,\alpha_k^{\vee}\rangle = h$, so that $\mu_k=s_{\alpha_k}(\lambda_k)$.
We now apply Proposition \ref{prop3.1} (with $\nu := \mu_k$, $\eta:= \lambda_k$, 
$\beta:= \beta_{r-k}$).
We have $\beta_{r-k}\in\Delta$ and 
$\beta_{r-k}\neq \alpha_k$ as $\alpha_k\not\in\Delta$. Furthermore, $m=
\langle s_{r-k}\cdots s_1\mu, \beta_{r-k}^{\vee}\rangle = 
-\langle \mu,  s_1\cdots s_{r-k-1}(\beta_{r-k})^{\vee}\rangle = a_{r-k}$. In the same
way $n=\langle \lambda_k, \beta_{r-k}^{\vee} \rangle =-b_{r-k}$. So by Proposition
\ref{prop3.1}, if we set
$$Y_{k+1} = y_{\beta_{r-k}}^{a_{r-k}} Y_k y_{\beta_{r-k}}^{b_{r-k}},$$
then we have an embedding $M(\mu_{k+1})= M(s_{\beta_{r-k}}\mu_k)\emb M(s_{\beta_{r-k}}
\lambda_k)=M(\lambda_{k+1})$ by $v_{\mu_{k+1}}\mapsto Y_{k+1} \cdot v_{\lambda_{k+1}}$.\par
Finally we note that $\lambda_r=\lambda$, $\mu_r=\mu$.
\end{pf}

It is possible to reformulate the algorithm in such a way that it looks more
like the method from \cite{mff}.
The construction described in \cite{mff} works as follows. Write $s_{\alpha} = 
s_{\alpha_{i_1}}\cdots s_{\alpha_{i_t}}$, as a product of simple reflections. 
For $1\leq k\leq t$
set $m_k = \langle s_{\alpha_{i_{k+1}}}\cdots s_{\alpha_{i_t}}\lambda, 
\alpha_{i_k}^{\vee}\rangle$.
Then $Y=y_{i_1}^{m_1}\cdots y_{i_t}^{m_t}$ is an element of $U({\mf n}^-)$ and
we have an embedding $M(\mu)\emb M(\lambda)$ by $v_{\mu}\mapsto 
Y\cdot v_{\lambda}$. Now, using the same notation as in the
description of the algorithm, the expression we get is 
$$y_{\beta_1}^{a_1}\cdots y_{\beta_r}^{a_r}y_{\beta}^h y_{\beta_r}^{b_r}\cdots
y_{\beta_1}^{b_1}.$$
As remarked in the proof of Proposition \ref{prop3.2}, $s_{\alpha} = 
s_1\cdots s_r s_{\beta} s_r\cdots s_1$ (where again we write $s_i=s_{\beta_i}$).
Furthermore, $b_k = \langle
s_{k-1}\cdots s_1\lambda, \beta_k^{\vee}\rangle$, $h=\langle \lambda,\alpha^{\vee}
\rangle = \langle s_r\cdots s_1\lambda, \beta^{\vee}\rangle$, $a_k = 
\langle s_{k+1}\cdots s_rs_{\beta}s_r\cdots s_1 \lambda, \beta_k^{\vee}\rangle$.
So we see that our method is a special case of the 
construction in \cite{mff}. However, the difference is that we have an
explicit method to rewrite the element above to an element of $U({\mf n}^-)$. 
By the next lemma the expression we use for $s_{\alpha}$ is the shortest possible
(so we cannot do essentially better by taking a different reduced expression).

\begin{lem}
The expression
$s_{\alpha} = s_1\cdots s_r s_{\beta} s_r\cdots s_1$ obtained by the first
step of the algorithm, is reduced.
\end{lem}

\begin{pf}
Set $\gamma = s_1(\alpha) = \alpha - m \beta_1$, where $m>0$. Then 
$s_{\alpha} = s_{s_1(\gamma)} = s_1s_{\gamma}s_1$. By induction, the 
expression $s_{\gamma} = s_2\cdots s_r s_{\beta} s_r\cdots s_2$ is reduced.
We show that $\ell( s_{\alpha} ) = \ell( s_{\gamma} ) + 2$. For this we use the
fact that the length of an element $w$ of the Weyl group is equal to the
number of positive roots that are mapped to negative roots by $w$. Write
$\Phi^+ = A\cup \{\beta_1\}$, where $A=\Phi^+\setminus \{\beta_1\}$. 
There is a positive root $\delta_0 \in \Phi$ with $s_{\gamma}s_1 \delta_0 =\beta_1$.
Set $S = \{ \delta \in A \mid s_{\gamma}s_1
\delta <0 \}\cup \{\delta_0, \beta_1\}$. Then $s_{\alpha}$ maps 
all elements of $S$ to negative roots. Since $\langle\gamma,\beta_1^{\vee}\rangle=-m<0$,
also $\langle \beta_1,\gamma^{\vee}\rangle <0$, and hence $s_{\gamma}(\beta_1)>0$.
So all roots that are mapped to negative roots by $s_{\gamma}$ are in $A$. Therefore,
since $s_1$ permutes $A$, there are $\ell(\gamma)$ roots $\delta\in A$ with 
$s_{\gamma}s_1(\delta)<0$. We conclude that  the cardinality of $S$ is 
$\ell(\gamma)+2$. So $\ell(\alpha) \geq \ell(\gamma)+2$, but that
means that $\ell(\alpha) = \ell(\gamma)+2$.
\end{pf}
\par

We can use the algorithm described in this section to construct general formulas
for singular elements. More precisely, let $\gamma$ be a fixed root in the root
system of $\mf{g}$. Then by applying the formulas of Section \ref{sec2} symbolically
we can derive a formula that given arbitrary weights $\lambda, \mu$ such that
$\langle \lambda, \gamma^{\vee}\rangle \in \Z_{>0}$ and $\mu = s_{\gamma}(\lambda)$
produces an element $Y\in U(\mf{n}^-)_{\lambda-\mu}$ such that 
$Y\cdot v_{\lambda}$ is singular. We illustrate this with an example.

\begin{exa}
Suppose that $\mf{g}$ is of type $A_3$. We use the same basis of $\mf{n}^-$ as in
Example \ref{exa1}. We consider the root $\alpha+\beta+\gamma$. Let $\lambda =
(r_1,r_2,r_3)$ be such that $h=r_1+r_2+r_3$ is a positive integer. A reduced
expression of $s_{\alpha+\beta+\gamma}$ is $s_{\alpha}s_{\beta}s_{\gamma}s_{\beta}
s_{\alpha}$. The corresponding element of $U(\mf{n}^-)$ is
$$Y=y_{\alpha}^{r_2+r_3}y_{\beta}^{r_3}y_{\gamma}^h y_{\beta}^{r_1+r_2}y_{\alpha}^{r_1}.$$
First we have 
$$y_{\beta}^{r_3}y_{\gamma}^h y_{\beta}^{r_1+r_2} = \sum_{k=0}^h (-1)^k 
\binom{h}{k} \binom{r_1+r_2}{k}k! y_{\beta}^{h-k}y_{\gamma}^{h-k}y_{\beta+\gamma}^k.$$
Now to obtain the formula for $Y$ we have to apply Lemma \ref{lem2} three times 
(and Lemma \ref{lem1} a few times), to obtain
\begin{align*}
Y=\sum_{k=0}^h\sum_{l_0}^k \sum_{s=0}^{h-k}\sum_{t=0}^s & (-1)^{k+l+s}
\binom{h}{k}\binom{r_1+r_2}{k}\binom{k}{l}\binom{r_1}{l}\binom{h-k}{s}\binom{r_1-l}{s}
\binom{s}{t}\binom{h-k}{t} \\
& k!l!s!t! y_{\alpha}^{h-l-s} y_{\beta}^{h-k-s}
y_{\gamma}^{h-k-t} y_{\alpha+\beta}^{s-t} y_{\beta+\gamma}^{k-l}
y_{\alpha+\beta+\gamma}^{l+t}.
\end{align*}
\end{exa}
Table \ref{tab2} contains a few running times of the implementation of this algorithm
in {\sf GAP}4.
\begin{table}[htb]
\begin{center}
\begin{tabular}{|c|r|r|}
\hline
type & length  & time (s) \\
\hline
$A_6$ & $29$ & 0.2 \\
$D_6$ & $109$ & 1.6\\
$E_6$ & $316$ & 2.9\\
$E_7$ & $2866$ & 26.3\\
$E_8$ & $>10556$ & $\infty$\\
\hline
\end{tabular}\caption{Running times for the computation of a formula for a
singular vector. 
}\label{tab2}
\end{center}
\end{table}
The root $\gamma$ is in each case the highest root of the root system. The length
of a formula is the number of summations it contains (so the length of the above
formula for $A_3$ is $4$). The computation for $E_8$ did not terminate in the available
amount of memory (100M). When the program exceeded the memory, the expression contained
$10556$ summations.

\begin{rem}
It is also possible to use this method to obtain formulas for a fixed type, but variable
rank. However, for that a convenient Chevalley basis needs to be chosen. We refer to
\cite{mff}, Section 5, for the formula for $A_n$. 
\end{rem}

\begin{rem}
We have chosen $\Q$ as the ground field, because it is easy to work with. However,
from the algorithm it is clear that instead we can choose any field $F$
of characteristic zero and construct embeddings of Verma modules with highest
weights from $FP$.
\end{rem}

\section{Composition of embeddings}\label{sec4}

In this section we consider the problem of obtaining an embedding $M(\nu)\emb M(\lambda)$,
where $\nu = s_{\alpha}s_{\beta}(\lambda) < s_{\beta}(\lambda) <\lambda$. The obvious 
way of doing this is to set $\mu = s_{\beta}(\lambda)$ and 
obtain the embeddings $M(\nu)\emb M(\mu)$, $M(\mu)\emb M(\lambda)$ 
and composing them. This amounts
to multiplying two elements of $U({\mf n}^-)$. This then corresponds to an
expression for $s_{\alpha}s_{\beta}$, which is not necessarily reduced. The
question arises whether in this case it is possible to do better, i.e., to
start with a reduced expression for $s_{\alpha}s_{\beta}=s_{\alpha_{i_1}}\cdots
s_{\alpha_{i_r}}$, set $m_k = \langle s_{\alpha_{i_{k+1}}}\cdots s_{\alpha_{i_t}}\lambda, 
\alpha_{i_k}^{\vee}\rangle$, and rewrite $Y=y_{i_1}^{m_1}\cdots y_{i_t}^{m_t}$
to an element of $U(\mf{n}^-)$. The next example  shows that this does not always
work.

\begin{exa}
Let $\Phi$ be of type $F_4$, with simple roots $\alpha_1,\ldots,\alpha_4$ and
Cartan matrix 
$$\begin{pmatrix} 2&-1&0&0\\ -1&2&-2&0\\ 0&-1&2&-1\\ 0&0&-1&2\end{pmatrix}.$$
Let $\alpha = \alpha_1+\alpha_2+2\alpha_3$ and $\beta=\alpha_1+2\alpha_2
+2\alpha_3+\alpha_4$. We abbreviate a weight $a_1\lambda_1+\cdots +a_4\lambda_4$
by $(a_1,a_2,a_3,a_4)$.
Set $\lambda = (\frac{5}{6},-\frac{1}{2},\frac{2}{3},0)$ and
$\nu = (-\frac{1}{6},-\frac{1}{2},-\frac{1}{3},2)$. Then $\nu = 
s_{\alpha}s_{\beta}(\lambda)$. Write $s_i=s_{\alpha_i}$. Then a reduced 
expression of $s_{\alpha}s_{\beta}$ is 
$$s_1s_2s_1s_3s_2s_1s_3s_2s_4s_3s_2s_1s_3s_2.$$ 
We get 
$$Y=y_1^{\frac{1}{6}} y_2^{\frac{2}{3}} y_1^{\frac{1}{2}} y_3^{\frac{5}{3}}  
y_2^{\frac{3}{2}} y_1 y_3^{\frac{4}{3}} y_2^{\frac{5}{6}} y_4 y_3^{\frac{4}{3}}
y_2^{\frac{1}{2}} y_1^{\frac{1}{3}} y_3^{-\frac{1}{3}} y_2^{-\frac{1}{2}}.$$
And I do not see any direct way to rewrite this as an element of $U({\mf{n}^-})$.\par
In general we have to obtain the embedding by composition. In this example
set $\mu= s_{\beta}(\lambda) = \lambda - \beta = (\frac{5}{6}, -\frac{3}{2},
\frac{5}{3},0)$. Then for the embedding $M(\mu)\emb M(\lambda)$ we get 
$$Y_1 = y_2^{\frac{2}{3}}y_3^{\frac{4}{3}}y_1^{\frac{2}{3}}y_2^{\frac{1}{2}}
y_3^{-\frac{1}{3}} y_4 y_3^{\frac{4}{3}} y_2^{\frac{1}{2}}y_1^{\frac{1}{3}}
y_3^{-\frac{1}{3}} y_2^{-\frac{1}{2}}.$$
For the embedding $M(\nu)\emb M(\mu)$ we get 
$$Y_2 = y_1^{\frac{1}{6}} y_3^{\frac{1}{3}} y_2 y_3^{\frac{5}{3}} y_1^{\frac{5}{6}}.$$
Then the product $Y_2Y_1$ will provide the embedding $M(\nu)\emb M(\lambda)$.
\end{exa} 

\section{Affine algebras}\label{sec5}

In this section we comment on finding embeddings of Verma modules of affine
Kac-Moody algebras. First we fix some notation and recall some facts. Our main
reference for this is \cite{kac}.\par
We let $\hat{\mf{g}}$ be the (untwisted) affine Lie algebra corresponding 
to $\mf{g}$, i.e.,
$$ \hat{\mf{g}} = \Q[t,t^{-1}] \otimes \mathfrak{g} \oplus \Q K \oplus \Q d$$
with multiplication
$$[ t^m\otimes x + a_1 K + b_1 d, t^n\otimes y + a_2 K + b_2 d ] =
(t^{m+n}\otimes [x,y] + b_1 nt^n \otimes y - b_2 mt^m \otimes x ) + 
m\delta_{m,-n} \kappa(x,y) K,$$
where $m,n\in \Z$, $x,y\in \mathfrak{g}$, $a_1,a_2,b_1,b_2\in \Q$ and 
$\kappa(~,~)$ is the Killing form on $\mathfrak{g}$. \par
The Lie algebra $\hat{\mf{g}}$ has a triangular decomposition $\hat{\mf{g}} = 
\hmf{n}^- \oplus \hmf{h} \oplus \hmf{n}^+$. Here $\hmf{n}^-$ is spanned by the 
$t^m \otimes y_i$ for $m\leq 0$, along with $t^n \otimes x_i$, and $t^n\otimes h_j$ 
for $n<0$. The subalgebra $\hmf{h}$ is spanned by the $t^0\otimes h_i$ and $K$ and $d$. 
Furthermore, $\hmf{n}^+$ is spanned by the $t^m\otimes x_i$ for $m\geq 0$, along with 
$t^n\otimes y_i$ and $t^n\otimes h_j$ for $n>0$. \par
The Verma module $M(\lambda)$ of highest weight $\lambda$ is defined in the same 
way as for $\mathfrak{g}$. As vector spaces $M(\lambda)\cong U(\hmf{n}^-)$.
Let $\alpha$ be a positive root of $\hmf{g}$. Then from \cite{kackazhdan} we 
get that $M(\lambda-n\alpha)$ embeds in $M(\lambda)$ if and only if  
$2(\lambda,\alpha) = n(\alpha,\alpha)$, where $n$ is  a positive integer.\par
Now if $\alpha$ is a real root with $2(\lambda,\alpha) = n(\alpha,\alpha)$,
then we can construct a singular vector in $U(\hmf{n}^-)_{n\alpha}$ by 
essentially the same method as in Section \ref{sec3}. The only difference
is the algorithm for rewriting $f_i^{n-r}af_i^r$, where $r\in \Q$, $a\in U(\hmf{n}^-)$,
and $f_i$ a basis element of $\hmf{n}^-$.
We need commutation relations $y^m f_i^r = f_i^r y^m +\cdots$, where $y$ 
runs through the basis elements of $\hmf{n}^-$. \par
First of all, if $f_i = t^j\otimes x_{\alpha}$ for some $\alpha\in\Phi$,
and $y=t^k\otimes x_{\beta}$ for some $\beta\in \Phi$ such that $\alpha+\beta
\in \Phi$, then set $y_{m\alpha+n\beta} = t^{mj+nk}\otimes x_{m\alpha+n\beta}$.
Set $B = \{ y_{m\alpha+n\beta} \mid m\alpha + n\beta \in \Phi \}$.
Then $B$ spans a subalgebra of $\hmf{g}$ isomorphic to the subalgebra
of $\mathfrak{g}$ spanned by the corresponding $x_{m\alpha+n\beta}$.
The isomorphism is given by $y_{m\alpha+n\beta}\mapsto x_{m\alpha+n\beta}$.
So we get the same formula as in the finite-dimensional case.\par
Now suppose that $\alpha+\beta=0$. Then $j+k\leq 0$; so $[f_i,y] = 
t^{j+k}\otimes h_{\alpha}$, where  $h_{\alpha} = [x_{\alpha},x_{-\alpha}]$. 
In this case we use the following relation:
\begin{align*}
(t^k\otimes x_{-\alpha})(t^j\otimes x_{\alpha})^r = &
(t^j\otimes x_{\alpha})^r(t^k\otimes x_{-\alpha}) - \\
& r (t^j\otimes x_{\alpha})^{r-1}(t^{k+j}\otimes h_{\alpha}) - r(r-1)
(t^j\otimes x_{\alpha})^{r-2}(t^{k+2j}\otimes x_{\alpha}),
\end{align*}
which is easily proved by induction. If $t^k\otimes x_{-\alpha}$ occurs with an exponent
$>1$ then we use this formula repeatedly. \par
The last possibility is 
$$(t^k\otimes h_q)(t^j\otimes x_{\alpha})^r = (t^j\otimes x_{\alpha})^r
(t^k\otimes h_q) +
r\langle \alpha,\alpha_q^{\vee}\rangle (t^j\otimes x_{\alpha})^{r-1}
(t^{k+j}\otimes x_{\alpha}).$$
Again, we use this formula repeatedly if $t^k\otimes h_q$ occurs with exponent $>1$.\par
Now we suppose that $\alpha = m\delta$ is an imaginary root with $(\lambda,\alpha)=0$
(here $\delta$ is the fundamental imaginary root).
Then $M(\lambda-n\alpha)\emb M(\lambda)$ for all positive integers $n$.
In this case there are a lot of singular elements. One class of them is easily 
constructed. Let $u_1,\ldots, u_q$, $u^1,\ldots, u^q$ be two basis of $\mathfrak{g}$,
dual to each other with respect to the Killing form. For $n>0$ set 
$$S_n = \sum_{i=1}^q \sum_{j=0}^n (t^{-j}\otimes u_i)( t^{j-n}\otimes u^i).$$

\begin{lem}\label{lemdelta}
Suppose that $(\lambda,\delta)=0$, then
$S_n\cdot v_{\lambda}$ is a singular vector of weight $n\delta$ in $M(\lambda)$.
\end{lem}

\begin{pf}
From \cite{kac}, 12.8 we have the Sugawara operators
$$T_s = \sum_{m\in \Z} \sum_{i=1}^q (t^{-m}\otimes u_i)(t^{m+s}\otimes u^i).$$
It is straightforward to see that $S_n\cdot v_{\lambda} = T_{-n}\cdot v_{\lambda}$.
Now $K$ acts on $M(\lambda)$ as scalar multiplication by 
$-h^{\vee}$, where $h^{\vee}$ is the dual Coxeter number. 
But also by \cite{kac}, Lemma 12.8 we have for $x\in \mathfrak{g}$:
$$[ t^m\otimes x , T_{-n}] = 2m(K+h^{\vee})(t^{m-n}\otimes x).$$
From this it follows that $x \cdot S_n v_{\lambda}=0$ for $0\leq i\leq l$, 
$x\in \mf{n}^+$. Therefore $S_n\cdot v_{\lambda}$ is a singular vector.
\end{pf}

Lemma \ref{lemdelta} provides an infinite number of singular vectors.
However, these are not the only ones. In \cite{malikov} it is shown that
for $n>0$ and $1\leq i\leq l$ there are independent elements 
$S_n^i\in U(\hmf{n}^-)$ of weight $n\delta$, such that $S_n^i\cdot v_{\lambda}$ is a 
singular vector. These $S_n^i$ are constructed from the generators of the 
centre of $U(\mf{g})$. In this construction the $S_n$ correspond to the
Casimir operator. However, with the exception of the Casimir operator, I
do not know of efficient algorithms to construct the generators of the centre of 
$U(\mf{g})$. For example, the explicit expressions  given in \cite{gauger}
for a generator of the centre of degree $s$ involve sums of $(\dim \mf{g})^s$ terms.
So constructing generators of the centre of $U(\mf{g})$
appears to be a very hard algorithmic problem in its own right. \par
The conclusion is that we have efficient algorithms to construct an inclusion
$M(\lambda-n\alpha)\emb M(\lambda)$ if $\alpha$ is a real root, or when 
$\alpha$ is imaginary. However, in the last case there are many singular 
vectors for which at present we have no efficient way of constructing them. 

\section{Constructing irreducible representations}\label{sec6}

In \cite{litt7}, P. Littelmann proves a theorem describing a particular
basis of the irreducible representations of $\mf{g}$, using
inclusions of Verma modules. Apart from giving a basis this result
also allows one to construct the irreducible representations of $\mf{g}$.
In this section we first briefly indicate how this works, and then give some 
experimental data concerning this algorithm. \par
The first ingredient of the construction is Littelmann's path method.
Here we only give a very rough description of that method; for the 
details we refer to \cite{litt6}, \cite{litt2}. A path is a piecewise 
linear function $\pi : [0,1] \to \R P$, such that $\pi(0)=0$.
Such a path is given by two
sequences $\bar{\mu}=(\mu_1,\ldots, \mu_r)$ and $\bar{a} = (a_0=0,a_1,
\ldots, a_r=1)$, where the $\mu_i\in \R P$ and the $a_i$ are real numbers
with $0=a_0<a_1<\ldots <a_r=1$. The path $\pi$ corresponding to this data
is given by 
$$\pi(t) = (t-a_{s-1})\mu_s + \sum_{i=1}^{s-1}(a_i-a_{i-1})\mu_i
\text{ ~ for } a_{s-1}\leq t\leq a_s.$$
Let $\lambda$ be a dominant weight. Then the path $\pi_{\lambda}$
is given by the sequences $(\lambda)$ and $(0,1)$, i.e., it is the straight
line from the origin to $\lambda$. For $\alpha\in \Delta$ there is a
path-operator $f_{\alpha}$. Given a path $\pi$, $f_{\alpha}(\pi)$ is a new path,
or $0$. Set $B(\lambda) = \{ f_{\alpha_{i_1}}\cdots f_{\alpha_{i_k}}
(\pi_{\lambda}) \mid k\geq 0, ~ \alpha_{i_j}\in \Delta\}$, and let
$V(\lambda)$ be the irreducible $\mf{g}$-module with highest weight
$\lambda$. Then from \cite{litt6}, \cite{litt2} we have that 
the endpoints of the paths in $B(\lambda)$ are weights of $V(\lambda)$
and the number of paths with endpoint $\mu$ is equal to the dimension
of the weight space in $V(\lambda)$ with weight $\mu$. \par
Let $\pi \in B(\lambda)$ be given by the sequences $(\mu_1,\ldots, 
\mu_r)$ and $(a_0=0,a_1,\ldots, a_r=1)$. Set $\mu_{r+1}=\lambda$ and
$\nu_i = a_i \mu_i$ and $\eta_i = a_i\mu_{i+1}$ for $1\leq i\leq r$. Then
it can be shown that $M(\nu_i)\emb M(\eta_i)$. Let $\Theta_i \in 
U(\mf{n}^-)_{\eta_i-\nu_i}$ be such that $\Theta_i \cdot v_{\eta_i}$ is a singular 
vector. Then set $\Theta_{\pi} = \Theta_1\cdots \Theta_r$. The element
$\Theta_{\pi}\in U(\mf{n}^-)_{\lambda-\pi(1)}$ is determined upto a 
multiplicative constant. \par
Now in \cite{litt7} an inclusion $B( m\lambda)\emb B(n\lambda)$ is
described for $m<n$. With this inclusion we can view $B(m\lambda)$ as 
as a subset of $B(n\lambda)$. Furthermore, $B(\lambda,\infty)$ denotes the
union of all $B(m\lambda)$ for $m\geq 1$. Write $\lambda = n_1\lambda_1+
\cdots +n_l\lambda_l$, and let $I(\lambda)$ be the left ideal of
$U(\mf{n}^-)$ generated by the elements $y_i^{n_i+1}$ for $1\leq i\leq l$.
Then $V(\lambda) = M(\lambda+\rho)/I(\lambda)\cdot v_{\lambda}$, where
$\rho = \lambda_1+\cdots +\lambda_l$. Now from \cite{litt7} we have the 
following result.

\begin{prop}
Suppose that all $n_i>0$. Then the set $\{\Theta_{\pi} \mid \pi \in 
B(\lambda,\infty), \pi\not\in B(\lambda) \}$ is a basis of $I(\lambda)$.
\end{prop}

(If some $n_i=0$ then there is a similar result, which we will omit here, cf.
\cite{litt7}.) \par
In order to construct and work with the quotient $M(\lambda+\rho)/I(\lambda)$,
we need a basis of $I(\lambda)$. If $\lambda-\mu$ is not a weight
of $V(\lambda)$, then $I(\lambda) \cap U(\mf{n}^-)_{\mu} =  U(\mf{n}^-)_{\mu}$.
So we only need bases of the spaces $I(\lambda) \cap U(\mf{n}^-)_{\mu}$ where
$\lambda-\mu$ is a weight of $V(\lambda)$. By the above theorem we can compute 
those bases by first computing paths $\pi \in B(\lambda, \infty)$ with $\pi(1)=
\lambda-\mu$, and then constructing the corresponding $\Theta_{\pi}$. We call this
algorithm A. \par
In Table \ref{tab1}, the running times are given of algorithm A on some
sample inputs. Also listed are the running times of the algorithm described
in \cite{gra5}, which uses a Gr\"{o}bner basis method to compute bases of
the spaces $I(\lambda) \cap U(\mf{n}^-)_{\mu}$. We call it algorithm B. 
In order to fairly compare both algorithms, the output in both cases consisted
of the representing matrices of a Chevalley basis of $\mf{g}$.

\begin{table}[htb]
\begin{center}
\begin{tabular}{|c|r|r|r|r|r|}
\hline
type & $\lambda$ & $\dim V(\lambda)$ & $\sharp$ inclusions & time A (s) & time B (s)\\
\hline
$A_2$ & $(2,2)$ & 27 & 64 & 1.0 & 1.3\\
$A_2$ & $(3,4)$ & 90 & 296 & 2.3 & 5.0\\
$A_2$ & $(5,5)$ & 216 & 788 & 6.5 & 14.2\\
\hline
$A_3$ & $(1,1,1)$ & 64 & 897 & 16.6 & 6.0\\
$A_3$ & $(2,1,1)$ & 140 & 2834 & 56.4 & 15.0\\
$A_3$ & $(2,1,2)$ & 300 & 7837 & 178.4 & 40.0\\
\hline
$B_2$ & $(2,2)$ & 81 & 807 & 10 & 6\\
$B_2$ & $(3,3)$ & 256 & 3330 & 56 & 23\\
$B_2$ & $(4,4)$ & 625 & 9502 & 347 & 79\\
\hline
\end{tabular}\caption{Running times (in seconds) of the algorithms A and B 
for the construction of $V(\lambda)$. The fourth column displays the number
of inclusions of Verma modules computed by algorithm A. The ordering of the 
fundamental weights is as in \cite{bou4}.}\label{tab1}
\end{center}
\end{table}

We see that for type $A_2$, algorithm A competes well with algorithm B. However,
for the other types considered this is not the case. In these cases huge numbers
of inclusions of Verma modules have to be constructed, which slows the algorithm
down considerably. I have also tried to construct $V(\lambda)$ for $\lambda =
(1,1,1)$, and $\mf{g}$ of type $B_3$. But algorithm A did not complete this 
calculation within the available amount of memory (100M).

\bibliographystyle{plain}

\def\cprime{$'$} \def\cprime{$'$} \def\cprime{$'$}

\end{document}